\documentclass[oneside,10pt]{article}
\usepackage{amsmath}
\usepackage{amssymb}
\usepackage{amsthm}

\newtheorem{theorem}{Theorem}[section]

\theoremstyle{definition}
\newtheorem{definition}{Definition}
\newtheorem{example}{Example}

\begin{document}

\begin{center}
\emph{Some examples on automorphisms of structural matrix algebras\bigskip }

\emph{Emira AKKURT and Mustafa Akkurt\bigskip }

\emph{Department of Mathematics, Gebze Institute of Technology, Kocaeli,
Turkey\bigskip }
\end{center}

Abstract: Automorphisms of structural matrix algebras in block upper
triangular form has been studied recently in \cite{Akkurt E-M Barker 2}, and
this work is a follow-up paper of that study. The aim of this paper is to
explain the topic in a much clearer and more comprehensible way by
presenting some examples.\emph{\bigskip }

\textbf{Keywords}: structural matrix algebras, block upper triangular\
matrix, automorphisms.\emph{\bigskip }


\section{\protect\Large Introduction and notation}

The topic of structural matrix algebras and automorphisms have been studied
in many papers, especially in \cite{Akkurt E-M Barker 1, Akkurt E-M Barker
2, Akkurt, Coelho}. The recent study \cite{Akkurt E-M Barker 2} gives an
answer to the open question given in \cite{Barker1} in 1993. For upper
triangular matrices we know that an automorphism is a composition of an
inner automorphism with a permutation similarity (see \cite{Barker1, Barker2}%
). Since structural matrix algebra is isomorphic with a block upper
triangular matrix algebra it is shown in \cite{Akkurt E-M Barker 2} that any
automorphism is a composition of an inner automorphism, a permutation
similarity and an automorphism generated by a transitive function $g$ on $%
\rho $ using block upper triangular form of structural matrix algebra. Here
we would like to give some application problems to express the main results
of \cite{Akkurt E-M Barker 2} in a numerical way. This study also can be
compared by \cite{Coelho}. As we mention in \cite{Akkurt E-M Barker 2}, we
use the block triangular form while Coelho uses graph theory in her paper.
Therefore we would like to use also same relations $\rho $ in her two
examples in \cite{Coelho} to show getting same results by using our approach
and technique.

Let $M_{n}(F)$ be the algebra of all $n\times n$ matrices over a field $F$
and let $(\{1,\ldots ,n\},\rho )$ be a quasi-ordered set (i.e. $\rho $ is
reflexive and transitive relation on the set $\{1,\ldots ,n\}$). The set 
\begin{equation*}
M_{n}(F,\rho )=\{A\in M_{n}(F):a_{ij}=0\text{ whenever }\left( i,j\right)
\notin \rho \}
\end{equation*}%
is subalgebra of $M_{n}(F)$ and we call $M_{n}(F,\rho )$ the algebra of $%
n\times n$ structural matrices over $F$ (with identity $I$).

Before working on examples, we would like to remark briefly to the
construction of block upper triangular form of structural matrix algebras.
We also need to remind some important results and main theorems given in 
\cite{Akkurt E-M Barker 2}. $M_{n}=M_{n}\left( F,\rho \right) $ which we
take to be in block upper triangular form that is constructed in \cite%
{Akkurt E-M Barker 2} and an automorphism $\Phi \in $Aut$\left( M_{n}\right)
.$ Let $F^{\ast }=F\diagdown \left\{ 0\right\} .$

\begin{definition}
A function $g:\rho \rightarrow F^{\ast }$ is transitive if and only if%
\begin{equation*}
g\left( i,j\right) g\left( j,k\right) =g\left( i,k\right)
\end{equation*}%
for all $\left( i,j\right) ,\left( j,k\right) \in \rho .$
\end{definition}

Every transitive function $g:\rho \rightarrow F^{\ast }$ determines an
automorphism $G\in $Aut$\left( M_{n}\left( F,\rho \right) \right) $ by
defining 
\begin{equation*}
G\left( E^{ij}\right) =g\left( i,j\right) E^{ij},\left( i,j\right) \in \rho ,
\end{equation*}%
where as we recall $E^{ij}$ is the $n\times n$ matrix with a $1$ in position 
$(i,j)$ and zeros elsewhere.

\section{{\protect\Large Construction of} {\protect\Large block upper
triangular form of structural matrix algebras}}

For a structural matrix algebra $M_{n}\left( F,\rho \right) $ we define an
equivalence relation $\overline{\rho }$ by 
\begin{equation*}
\left( i,j\right) \in \overline{\rho }\text{ if and only if }\left(
i,j\right) ,\left( j,i\right) \in \rho .
\end{equation*}%
Let $\left[ r_{1}\right] ,\left[ r_{2}\right] ,\ldots ,\left[ r_{p}\right] $
denote the distinct equivalence classes of $\overline{\rho }$ with
representatives $r_{1},r_{2},\ldots ,r_{p}.$ And we construct a permutation $%
\pi $ in \cite{Akkurt E-M Barker 2} as follows. For 
\begin{equation*}
\left[ r_{1}\right] =\left\{ r_{1}=r_{11},r_{12},\ldots ,r_{1m_{1}}\right\}
\end{equation*}%
let 
\begin{equation*}
\pi \left( 1\right) =\pi \left( r_{11}\right) =1,\pi \left( r_{12}\right)
=2,\ldots ,\pi \left( r_{1m_{1}}\right) =m_{1}
\end{equation*}%
and in general if 
\begin{equation*}
\left[ r_{k}\right] =\left\{ r_{k1},r_{k2},\ldots ,r_{km_{k}}\right\}
\end{equation*}%
then 
\begin{equation*}
\pi \left( r_{kj}\right) =m_{1}+m_{2}+\cdots +m_{k-1}+j.
\end{equation*}%
If we apply this permutation similarity to $M_{n}\left( F,\rho \right) $ we
have the following relation $\rho ^{\prime }$ where 
\begin{equation*}
\left( i,j\right) \in \rho ^{\prime }\iff \left( \pi ^{-1}\left( i\right)
,\pi ^{-1}\left( j\right) \right) \in \rho
\end{equation*}%
and $M_{n}\left( F,\rho ^{\prime }\right) $ consists of all matrices of the
form%
\begin{equation*}
\begin{bmatrix}
M_{m_{1}}\left( F\right) & M_{m_{1}\times m_{2}}\left( R\right) & \cdots & 
M_{m_{1}\times m_{p}}\left( R\right) \\ 
0 & M_{m_{2}}\left( F\right) & \cdots & M_{m_{2}\times m_{p}}\left( R\right)
\\ 
\vdots & 0 & \ddots & \vdots \\ 
0 & \cdots & 0 & M_{m_{p}}\left( F\right)%
\end{bmatrix}%
,
\end{equation*}%
where $R$ is either $F$ or $0.$

And we define%
\begin{equation}
\left[ i\right] _{\rho ^{\prime }}\leq \left[ j\right] _{\rho ^{\prime }} 
\tag{1}  \label{star}
\end{equation}%
to mean $\left( i,j\right) \in \rho ^{\prime }$. Then $\forall p_{1}\in %
\left[ i\right] ,p_{2}\in \left[ j\right] ,$%
\begin{equation*}
\left( p_{1},p_{2}\right) \in \rho ^{\prime }\text{ if and only if }\left(
i,j\right) \in \rho ^{\prime }.
\end{equation*}

Next, let $\left[ t_{1}\right] ,\left[ t_{2}\right] ,\ldots ,\left[ t_{q}%
\right] $ be the classes which are incomparable with any other class, i.e.,
for any $t_{k},$ $k=1,2,\ldots ,q,\ $for any $r_{j}\neq t_{k}$ we have 
\begin{eqnarray*}
\left[ r_{j}\right] &\nleqslant &\left[ t_{k}\right] \ \ \text{ }\left( 
\text{that is }\left( r_{j},t_{k}\right) \notin \rho ^{\prime }\right) \text{
and} \\
\left[ t_{k}\right] &\nleqslant &\left[ r_{j}\right] \ \ \ \left( \text{that
is }\left( t_{k},r_{j}\right) \notin \rho ^{\prime }\right) .
\end{eqnarray*}%
Then we performed a permutation similarity of the block triangular matrix
algebras which makes the diagonal blocks indexed by the indices of $\left[
t_{1}\right] ,\left[ t_{2}\right] ,\ldots ,\left[ t_{q}\right] $ the last $q$
diagonal blocks. Denote the corresponding permutation by $\pi _{1}$ and the
quasi-order relation by $\sigma .$ Further the last $q$ $\sigma $-classes,
which we denote by $\left[ t_{1}\right] _{\sigma },\left[ t_{2}\right]
_{\sigma },\ldots ,\left[ t_{q}\right] _{\sigma }$, are such that no one of
these classes is comparable to any other $\sigma $-class (with $\leq $
defined by (\ref{star}) and $\rho ^{\prime }$ replaced by $\sigma $). This
means that $R=0$ for all the $M_{m_{j}\times m_{k}}\left( R\right) $ with $%
k=p-q+1,\ldots ,p.$

To summarize this, the classes comparable to some other class are $\left[
r_{j}\right] $ $\left( j=1,\ldots ,l\right) $ while the classes comparable
to no other class are $\left[ t_{k}\right] $ $\left( k=1,\ldots ,q\right) $
each of cardinality $m_{l}$ and $m_{q}$ respectively. Also $l+q=p.$ Hence we
have our diagonal form 
\begin{equation*}
\begin{bmatrix}
M_{m_{1}}(F) & M_{m_{1}\times m_{2}}(R) & \cdots & M_{m_{1}\times m_{l}}(R)
& 0 & \cdots & 0 \\ 
0 & M_{m_{2}}(F) & \cdots & M_{m_{2}\times m_{l}}(R) & 0 & \cdots & 0 \\ 
\vdots & \ddots & \ddots & \vdots & \vdots & \vdots & \vdots \\ 
0 & 0 & \cdots & M_{m_{l}}(F) & 0 & \cdots & 0 \\ 
0 & 0 & 0 & 0 & M_{m_{l+1}}(F) &  & 0 \\ 
\vdots & \vdots & \vdots & \vdots & \vdots & \ddots & \vdots \\ 
0 & 0 & \cdots & 0 & 0 & \cdots & M_{m_{p}}(F)%
\end{bmatrix}%
.
\end{equation*}%
Note that the way we have chosen the $r_{j}$ and $t_{k}$ in the construction
of $\pi _{1}$ we have that


\begin{enumerate}
\item[(i)] each $r_{j}$ and each $t_{k}$ is minimal in its class,

\item[(ii)] \label{ii}$r_{1}<\ldots <r_{l}$ and without loss of generality
we can have $t_{1}<\ldots <t_{q}.$
\end{enumerate}


Now, if we consider a block $M_{m_{a}\times m_{b}}(R)$ with $a\leq l$ and $%
b\leq l,$ then the row and the column indices are $\pi _{1}$ classes, say $%
\left[ r_{a}\right] $ and $\left[ r_{b}\right] $ respectively with $a<b.$ If 
$\left[ r_{a}\right] \leq \left[ r_{b}\right] $ we have $R=F;$ otherwise $%
R=0.$

For notational simplicity, we may write $M_{i}$ or $M_{ij}$ for $M_{m_{i}}$
or $M_{m_{i}\times m_{j}}$ respectively$.$ The subset of $M_{n}\left( F,\rho
\right) $ which has elements of $M_{ij}$ $\left( \text{or }M_{i}\right) $ in
the $i^{\text{th}}$ block row and $j^{\text{th}}$ block column (or $i^{\text{%
th}}$ block diagonal) and zero elsewhere will be denoted by $\overline{M}%
_{ij}.$ Analogously we let $A_{ij}$ denote an $m_{i}\times m_{j}$ matrix in $%
M_{ij}$ while $\overline{A}_{ij}$ is the corresponding element of $\overline{%
M}_{ij}.$ If $i=j$ we write $M_{i}$, $A_{i},\overline{M}_{i},$ and $%
\overline{A}_{i},$ for the corresponding sets and matrices.

\begin{theorem}
\cite{Akkurt E-M Barker 2}If $M_{n}\left( F,\rho \right) $ is semisimple and
if $\Phi \in $Aut$\left( M_{n}\right) ,$ then we can write%
\begin{equation*}
\Phi =\Psi \circ P_{\tau },
\end{equation*}%
where $\Psi $ is an inner automorphism, and $P_{\tau }$ is a permutation
similarity which is in Aut$\left( M_{n}\right) .$
\end{theorem}

\begin{theorem}
\cite{Akkurt E-M Barker 2}\label{FT}(Factorization Theorem) If $\Phi \in $Aut%
$\left( M_{n}\right) $ then we can write%
\begin{equation*}
\Phi =\Psi _{A}\circ G\circ P_{\tau },
\end{equation*}%
where $\Psi _{A}$ is an inner automorphism induced by $A$, $G$ is an
automorphism defined by transitive function $g$ on $\rho $ and $P_{\tau }$
is a permutation similarity which is in Aut$\left( M_{n}\right) .$
\end{theorem}

\section{\protect\Large Examples}

\begin{example}
Let $N=\left\{ 1,2,3,4,5,6\right\} $ and
\end{example}

$\rho =\{(1,1),(1,5),(1,6),\left( 2,2\right) ,\left( 2,3\right) ,\left(
3,2\right) ,\left( 3,3\right) ,\left( 4,4\right) ,\left( 5,1\right) ,\left(
5,5\right) ,(5,6),$\linebreak $\left( 6,1\right) ,\left( 6,5\right) ,(6,6)\}$

$\overline{\rho }=\rho $

then the structural matrix algebra is

\begin{equation*}
M_{6}\left( F,\rho \right) =\left\{ 
\begin{bmatrix}
F & 0 & 0 & 0 & F & F \\ 
0 & F & F & 0 & 0 & 0 \\ 
0 & F & F & 0 & 0 & 0 \\ 
0 & 0 & 0 & F & 0 & 0 \\ 
F & 0 & 0 & 0 & F & F \\ 
F & 0 & 0 & 0 & F & F%
\end{bmatrix}%
:F\text{ is a field}\right\} .
\end{equation*}

We can construct equivalence classes and define permutation $\pi $ as
follows,

$[r_{1}]=\{r_{11},r_{12},r_{13}\},\
[r_{2}]=\{r_{21},r_{22}\},[r_{3}]=\{r_{31}\},$

where $r_{1}=1,r_{12}=5,r_{13}=6,\ r_{21}=2,r_{22}=3,r_{31}=4$

$[1]=\{r_{1}=1,r_{12}=5,r_{13}=6\},[2]=\{r_{21}=2,r_{22}=3\},[4]=\{r_{31}=4%
\}.$

$\pi \left( 1\right) =1,\pi \left( 5\right) =2,\pi \left( 6\right) =3,\pi
\left( 2\right) =4,\pi \left( 3\right) =5,\pi \left( 4\right) =6,$

i.e. $\pi =(24635)\Rightarrow \pi ^{-1}=\left( 25364\right) $

Using $\left( i,j\right) \in \rho ^{\prime }\iff \left( \pi ^{-1}\left(
i\right) ,\pi ^{-1}\left( j\right) \right) \in \rho $,

So we have

$\rho ^{\prime }=\{(1,1),(1,2),(1,3),\left( 2,1\right) ,\left( 2,2\right)
,\left( 2,3\right) ,\left( 3,1\right) ,\left( 3,2\right) ,\left( 3,3\right)
,\left( 4,4\right) ,\left( 4,5\right) ,\linebreak \left( 5,4\right)
,(5,5),\left( 6,6\right) \}$ and the structural matrix algebra related to $%
\rho ^{\prime }$ is

\begin{equation*}
M_{6}\left( F,\rho ^{\prime }\right) =\left\{ 
\begin{bmatrix}
F & F & F & 0 & 0 & 0 \\ 
F & F & F & 0 & 0 & 0 \\ 
F & F & F & 0 & 0 & 0 \\ 
0 & 0 & 0 & F & F & 0 \\ 
0 & 0 & 0 & F & F & 0 \\ 
0 & 0 & 0 & 0 & 0 & F%
\end{bmatrix}%
:F\text{ is a field}\right\} .
\end{equation*}

Then the structural matrix algebra $M_{6}\left( F,\rho \right) $ is
permutation similarity to the structural matrix algebra $M_{6}\left( F,\rho
^{\prime }\right) .$

Since $M_{6}\left( F,\rho ^{\prime }\right) $ in block diagonal form and
hence semisimple, then we have $\overline{M}_{1},\overline{M}_{2},$ and $%
\overline{M}_{3}$ are ideals which are simple and we can follow the proof of
Theorem 1 in \cite{Akkurt E-M Barker 2} and write

\begin{equation*}
M_{6}\left( F,\rho ^{\prime }\right) =\overline{M}_{1}\oplus \overline{M}%
_{2}\oplus \overline{M}_{3}\text{, where}
\end{equation*}%
\begin{equation*}
\overline{M}_{1}=\left\{ 
\begin{bmatrix}
F & F & F & 0 & 0 & 0 \\ 
F & F & F & 0 & 0 & 0 \\ 
F & F & F & 0 & 0 & 0 \\ 
0 & 0 & 0 & 0 & 0 & 0 \\ 
0 & 0 & 0 & 0 & 0 & 0 \\ 
0 & 0 & 0 & 0 & 0 & 0%
\end{bmatrix}%
:F\text{ is a field}\right\} =\left\{ 
\begin{bmatrix}
A_{3} & 0 \\ 
0 & 0%
\end{bmatrix}%
:A_{3}\in M_{3}\left( F\right) \right\} ,
\end{equation*}

\begin{equation*}
\overline{M}_{2}=\left\{ 
\begin{bmatrix}
0 & 0 & 0 & 0 & 0 & 0 \\ 
0 & 0 & 0 & 0 & 0 & 0 \\ 
0 & 0 & 0 & 0 & 0 & 0 \\ 
0 & 0 & 0 & F & F & 0 \\ 
0 & 0 & 0 & F & F & 0 \\ 
0 & 0 & 0 & 0 & 0 & 0%
\end{bmatrix}%
:F\text{ is a field}\right\} =\left\{ 
\begin{bmatrix}
0 & 0 & 0 \\ 
0 & A_{2} & 0 \\ 
0 & 0 & 0%
\end{bmatrix}%
:A_{2}\in M_{2}\left( F\right) \right\} ,
\end{equation*}

\begin{equation*}
\ \overline{M}_{3}=\left\{ 
\begin{bmatrix}
0 & 0 & 0 & 0 & 0 & 0 \\ 
0 & 0 & 0 & 0 & 0 & 0 \\ 
0 & 0 & 0 & 0 & 0 & 0 \\ 
0 & 0 & 0 & 0 & 0 & 0 \\ 
0 & 0 & 0 & 0 & 0 & 0 \\ 
0 & 0 & 0 & 0 & 0 & F%
\end{bmatrix}%
:F\text{ is a field}\right\} =\left\{ 
\begin{bmatrix}
0 & 0 \\ 
0 & A_{1}%
\end{bmatrix}%
:A_{1}\in F\right\} .
\end{equation*}

Since, $A_{1},A_{2},$ and $A_{3}$ are of different size we have 
\begin{equation*}
\Phi \left( \overline{M}_{i}\right) =\overline{M}_{f\left( i\right) }=%
\overline{M}_{i}.
\end{equation*}%
i.e. $f$ is identity bijection and extent $f$ to a permutation $\pi _{1}$ on 
$\left\{ 1,\ldots ,6\right\} $ by way of (ii)\ and the corresponding
permutation matrix , then ($P_{\pi _{1}}=I$) 
\begin{equation*}
\pi _{1}\left( A\right) =P_{\pi _{1}^{-1}}AP_{\pi _{1}}
\end{equation*}%
is an identity automorphism of $M_{n}\left( F,\rho ^{\prime }\right) $ while 
$\Psi =\Phi $ is an automorphism of $M_{n}\left( F,\rho ^{\prime }\right) $
which is inner.

\begin{example}
Let $N=\left\{ 1,2,3\right\} $ and \newline
\end{example}

$\rho =\left\{ \left( 1,1\right) ,\left( 1,2\right) ,\left( 2,2\right)
,\left( 3,2\right) ,\left( 3,3\right) \right\} $

$\overline{\rho }=\left\{ \left( 1,1\right) ,\left( 2,2\right) ,\left(
3,3\right) \right\} $

and the structural matrix algebra is

\begin{equation*}
M_{3}\left( F,\rho \right) =\left\{ 
\begin{bmatrix}
a_{11} & a_{12} & 0 \\ 
0 & a_{22} & 0 \\ 
0 & a_{32} & a_{33}%
\end{bmatrix}%
:a_{ij}\in F\right\} .
\end{equation*}

We can construct equivalence classes and define permutation $\pi $ as
follows,

$[r_{1}]=\{r_{11}\},\ [r_{2}]=\{r_{21}\},[r_{3}]=\{r_{31}\},$ where $%
r_{1}=1,\ r_{21}=2,r_{31}=3$

$[1]=\{r_{1}=1\},[2]=\{r_{21}=2\},[3]=\{r_{31}=3\}.$

$\pi \left( 1\right) =1,\pi \left( 2\right) =3,\pi \left( 3\right) =2,\ $%
i.e., $\pi =(23)=\pi ^{-1}$

Using $\left( i,j\right) \in \rho ^{\prime }\iff \left( \pi ^{-1}\left(
i\right) ,\pi ^{-1}\left( j\right) \right) \in \rho $, we have,

$\rho ^{\prime }=\{(1,1),(1,3),\left( 2,2\right) ,\left( 2,3\right) ,(3,3)\}$
and the structural matrix algebra related to $\rho ^{\prime }$ is

\begin{equation*}
M_{3}\left( F,\rho ^{\prime }\right) =\left\{ 
\begin{bmatrix}
a_{11} & 0 & a_{13} \\ 
0 & a_{22} & a_{23} \\ 
0 & 0 & a_{33}%
\end{bmatrix}%
:F\text{ is a field}\right\} .
\end{equation*}

Then, the structural matrix algebra $M_{3}\left( F,\rho \right) $ is
permutation similarity to the structural matrix algebra $M_{3}\left( F,\rho
^{\prime }\right) .$

We now consider\ the permutation $\pi _{1}=(12)$ so

\begin{equation*}
P_{\pi _{1}}=%
\begin{bmatrix}
0 & 1 & 0 \\ 
1 & 0 & 0 \\ 
0 & 0 & 1%
\end{bmatrix}%
,P_{\pi _{1}^{-1}=}%
\begin{bmatrix}
0 & 1 & 0 \\ 
1 & 0 & 0 \\ 
0 & 0 & 1%
\end{bmatrix}%
\text{, take }B=%
\begin{bmatrix}
a_{11} & 0 & a_{13} \\ 
0 & a_{22} & a_{23} \\ 
0 & 0 & a_{33}%
\end{bmatrix}%
\in M_{3}\left( F,\rho ^{\prime }\right)
\end{equation*}

\begin{equation*}
P_{\pi _{1}^{-1}}BP_{\pi _{1}}=%
\begin{bmatrix}
0 & 1 & 0 \\ 
1 & 0 & 0 \\ 
0 & 0 & 1%
\end{bmatrix}%
\begin{bmatrix}
a_{11} & 0 & a_{13} \\ 
0 & a_{22} & a_{23} \\ 
0 & 0 & a_{33}%
\end{bmatrix}%
\begin{bmatrix}
0 & 1 & 0 \\ 
1 & 0 & 0 \\ 
0 & 0 & 1%
\end{bmatrix}%
=%
\begin{bmatrix}
a_{22} & 0 & a_{23} \\ 
0 & a_{11} & a_{13} \\ 
0 & 0 & a_{33}%
\end{bmatrix}%
\in M_{3}\left( F,\rho ^{\prime }\right) .
\end{equation*}

Therefore $\pi _{1}$ is an automorphism of $\rho ^{\prime }.$ Since every
transitive function is trivial, by the lemma 4.10 of \cite[Theorem B]{Coelho}%
, we have $G$ is inner automorphism, i.e., $G=\Psi _{\left( A\right) }$ for
some $A\in M_{3}\left( F,\rho ^{\prime }\right) $.

Using Lemma 4.8 in \cite[Theorem B]{Coelho} we have

$g\left( i,i\right) =1,\ \ i=1,2,3$ and $g\left( i,j\right) =s\left(
i\right) s^{-1}\left( j\right) $ where $s:\left\{ 1,2,3\right\}
\longrightarrow F^{\ast }$ and $\pi _{1}=(12)$.

Let $G=g\left( \pi _{1}\left( i\right) ,\pi _{1}\left( j\right) \right)
E^{\pi _{1}\left( i\right) \pi _{1}\left( j\right) }$ then

$G\left( E^{11}\right) =g\left( \pi _{1}\left( 1\right) ,\pi _{1}\left(
1\right) \right) E^{\pi _{1}\left( 1\right) \pi _{1}\left( 1\right) }=E^{22}$

$G\left( E^{22}\right) =g\left( \pi _{1}\left( 2\right) ,\pi _{1}\left(
2\right) \right) E^{\pi _{1}\left( 2\right) \pi _{1}\left( 2\right) }=E^{11}$

$G\left( E^{33}\right) =g\left( \pi _{1}\left( 3\right) ,\pi _{1}\left(
3\right) \right) E^{\pi _{1}\left( 3\right) \pi _{1}\left( 3\right) }=E^{33}$

$G\left( E^{13}\right) =g\left( \pi _{1}\left( 1\right) ,\pi _{1}\left(
3\right) \right) E^{\pi _{1}\left( 1\right) \pi _{1}\left( 3\right)
}=g\left( 2,3\right) E^{23}=s\left( 2\right) s^{-1}\left( 3\right)
E^{23}=bE^{23},$ where $b\in F^{\ast }$

$G\left( E^{23}\right) =g\left( \pi _{1}\left( 2\right) ,\pi _{1}\left(
3\right) \right) E^{\pi _{1}\left( 2\right) \pi _{1}\left( 3\right)
}=g\left( 1,3\right) E^{13}=s\left( 1\right) s^{-1}\left( 3\right)
E^{13}=aE^{13},$ where $a\in F^{\ast }$

$A=%
\begin{bmatrix}
1 & 0 & a \\ 
0 & 1 & b \\ 
0 & 0 & 1%
\end{bmatrix}%
,A^{-1}=%
\begin{bmatrix}
1 & 0 & -a \\ 
0 & 1 & -b \\ 
0 & 0 & 1%
\end{bmatrix}%
$ then $G=\Psi _{A}$ and

\begin{eqnarray*}
\left( \Psi _{A}\circ P_{\pi _{1}}\right) \left( B\right) &=&\Psi _{A}\left(
P_{\pi _{1}^{-1}}BP_{\pi _{1}}\right) =A^{-1}\left( P_{\pi _{1}^{-1}}BP_{\pi
_{1}}\right) A \\
&=&%
\begin{bmatrix}
1 & 0 & -a \\ 
0 & 1 & -b \\ 
0 & 0 & 1%
\end{bmatrix}%
\begin{bmatrix}
a_{22} & 0 & a_{23} \\ 
0 & a_{11} & a_{13} \\ 
0 & 0 & a_{33}%
\end{bmatrix}%
\begin{bmatrix}
1 & 0 & a \\ 
0 & 1 & b \\ 
0 & 0 & 1%
\end{bmatrix}
\\
&=&%
\begin{bmatrix}
a_{22} & 0 & a_{23}+aa_{22}-aa_{33} \\ 
0 & a_{11} & a_{13}+ba_{11}-ba_{33} \\ 
0 & 0 & a_{33}%
\end{bmatrix}%
\in M_{3}\left( F,\rho ^{\prime }\right) .
\end{eqnarray*}

If we take $\Phi :M_{3}\left( F,\rho ^{\prime }\right) \longrightarrow
M_{3}\left( F,\rho ^{\prime }\right) $ by%
\begin{equation*}
\Phi \left( 
\begin{bmatrix}
a_{11} & 0 & a_{13} \\ 
0 & a_{22} & a_{23} \\ 
0 & 0 & a_{33}%
\end{bmatrix}%
\right) =%
\begin{bmatrix}
a_{22} & 0 & a_{23}+aa_{22}-aa_{33} \\ 
0 & a_{11} & a_{13}+ba_{11}-ba_{33} \\ 
0 & 0 & a_{33}%
\end{bmatrix}%
\in M_{3}\left( F,\rho ^{\prime }\right) ,
\end{equation*}%
then 
\begin{equation*}
\Phi =\Psi _{A}\circ P_{\pi _{1}}.
\end{equation*}

\begin{example}
Let $N=\left\{ 1,2,3,4,5,6\right\} $ and \newline
\end{example}

$\rho =\{(1,1),(1,2),(1,3),(1,4),\left( 2,2\right) ,\left( 2,3\right)
,\left( 3,2\right) ,\left( 3,3\right) ,\left( 4,4\right) ,\left( 5,2\right)
,\left( 5,3\right) ,(5,4),(5,5),$

$(5,6),\left( 6,2\right) ,\left( 6,3\right) ,(6,4),(6,5),(6,6)\}\bigskip $

$\overline{\rho }=\{(1,1),\left( 2,2\right) ,\left( 2,3\right) ,\left(
3,2\right) ,\left( 3,3\right) ,\left( 4,4\right)
,(5,5),(5,6),(6,6)\}\bigskip $

then the structural matrix algebra is

\begin{equation*}
M_{6}\left( F,\rho \right) =\left\{ 
\begin{bmatrix}
F & F & F & F & 0 & 0 \\ 
0 & F & F & 0 & 0 & 0 \\ 
0 & F & F & 0 & 0 & 0 \\ 
0 & 0 & 0 & F & 0 & 0 \\ 
0 & F & F & F & F & F \\ 
0 & F & F & F & F & F%
\end{bmatrix}%
:F\text{ is a field}\right\} .
\end{equation*}

We can construct equivalence classes and define permutation $\pi $ as
follows,

\begin{equation*}
\lbrack r_{1}]=\{r_{1}=r_{11}\},\
[r_{2}]=\{r_{21},r_{22}\},[r_{3}]=\{r_{31}\},[r_{4}]=\{r_{41},r_{42}\},
\end{equation*}

where 
\begin{equation*}
r_{1}=1,r_{21}=2,r_{22}=3,r_{31}=4,r_{41}=5,r_{42}=6
\end{equation*}

$[1]=\{r_{1}=1\},[2]=\{r_{21}=2,r_{22}=3\},[4]=\{r_{31}=4\},[5]=%
\{r_{41}=5,r_{42}=6\}$

$\pi \left( 1\right) =1,\pi \left( 5\right) =2,\pi \left( 6\right) =3,\pi
\left( 4\right) =4,\pi \left( 2\right) =5,\pi \left( 3\right) =6,$ i.e. $\pi
=(25)(36)=\pi ^{-1}$

Using $\left( i,j\right) \in \rho ^{\prime }\iff \left( \pi ^{-1}\left(
i\right) ,\pi ^{-1}\left( j\right) \right) \in \rho $, we have,\newline

$\rho ^{\prime }=\{(1,1),(1,4),(1,5),(1,6),\left( 2,2\right) .\left(
2,3\right) ,(2,4),(2,5),\left( 2,6\right) ,(3,2),(3,3),\left( 3,4\right)
,\left( 3,5\right) ,\newline
(3,6),\left( 4,4\right) ,(5,5),\left( 6,6\right) \left( 6,5\right) ,\left(
6,6\right) \},$

and the structural matrix algebra related to $\rho ^{\prime }$ is

\begin{equation*}
M_{6}\left( F,\rho ^{\prime }\right) =\left\{ 
\begin{bmatrix}
F & 0 & 0 & F & F & F \\ 
0 & F & F & F & F & F \\ 
0 & F & F & F & F & F \\ 
0 & 0 & 0 & F & 0 & 0 \\ 
0 & 0 & 0 & 0 & F & F \\ 
0 & 0 & 0 & 0 & F & F%
\end{bmatrix}%
:F\text{ is a field}\right\} .
\end{equation*}

Then, the structural matrix algebra $M_{6}\left( F,\rho \right) $ is
permutation similarity to the structural matrix algebra $M_{6}\left( F,\rho
^{\prime }\right) .$Then

\begin{equation*}
\overline{M}_{1}=\left\{ 
\begin{bmatrix}
F & 0 & 0 & 0 & 0 & 0 \\ 
0 & 0 & 0 & 0 & 0 & 0 \\ 
0 & 0 & 0 & 0 & 0 & 0 \\ 
0 & 0 & 0 & 0 & 0 & 0 \\ 
0 & 0 & 0 & 0 & 0 & 0 \\ 
0 & 0 & 0 & 0 & 0 & 0%
\end{bmatrix}%
:F\text{ is a field}\right\} =\left\{ 
\begin{bmatrix}
A_{1} & 0 \\ 
0 & 0%
\end{bmatrix}%
:A_{1}\in F\right\} ,
\end{equation*}%
\begin{equation*}
\overline{M}_{2}=\left\{ 
\begin{bmatrix}
0 & 0 & 0 & 0 & 0 & 0 \\ 
0 & F & F & 0 & 0 & 0 \\ 
0 & F & F & 0 & 0 & 0 \\ 
0 & 0 & 0 & 0 & 0 & 0 \\ 
0 & 0 & 0 & 0 & 0 & 0 \\ 
0 & 0 & 0 & 0 & 0 & 0%
\end{bmatrix}%
:F\text{ is a field}\right\} =\left\{ 
\begin{bmatrix}
0 & 0 & 0 \\ 
0 & A_{2} & 0 \\ 
0 & 0 & 0%
\end{bmatrix}%
:A_{2}\in M_{2}\right\} ,
\end{equation*}

\begin{equation*}
\overline{M}_{3}=\left\{ 
\begin{bmatrix}
0 & 0 & 0 & 0 & 0 & 0 \\ 
0 & 0 & 0 & 0 & 0 & 0 \\ 
0 & 0 & 0 & 0 & 0 & 0 \\ 
0 & 0 & 0 & F & 0 & 0 \\ 
0 & 0 & 0 & 0 & 0 & 0 \\ 
0 & 0 & 0 & 0 & 0 & 0%
\end{bmatrix}%
:F\text{ is a field}\right\} =\left\{ 
\begin{bmatrix}
0 & 0 & 0 \\ 
0 & A_{3} & 0 \\ 
0 & 0 & 0%
\end{bmatrix}%
:A_{3}\in F\right\} ,
\end{equation*}%
\begin{equation*}
\overline{M}_{4}=\left\{ 
\begin{bmatrix}
0 & 0 & 0 & 0 & 0 & 0 \\ 
0 & 0 & 0 & 0 & 0 & 0 \\ 
0 & 0 & 0 & 0 & 0 & 0 \\ 
0 & 0 & 0 & 0 & 0 & 0 \\ 
0 & 0 & 0 & 0 & F & F \\ 
0 & 0 & 0 & 0 & F & F%
\end{bmatrix}%
:F\text{ is a field}\right\} =\left\{ 
\begin{bmatrix}
0 & 0 \\ 
0 & A_{4}%
\end{bmatrix}%
;A_{4}\in M_{2}\right\} .
\end{equation*}

We have a permutation when and only when two blocks are of the same size,
then we may permute the blocks $A_{i}$ and $A_{j}.$ Otherwise not. Since, $%
A_{1}$ and $A_{3},$ $A_{2}$ and $A_{4}$ are of the same size we are going to
check whether these clases may permute.

Now lets consider the blocks $A_{1}$ and $A_{3},$ i.e.$,$the equivalent
classes $[1]$ and $[4],$ and consider the automorphism $\pi _{1}$ such that $%
\pi _{1}\left( 1\right) =4$ and $\pi _{1}\left( 4\right) =1$ then $\left(
\pi _{1}\left( 1\right) ,\pi _{1}\left( 4\right) \right) =\left( 4,1\right)
\notin \rho ^{\prime }.$ Therefore these blocks can not permute so the
equivalent classes are not comparable, i.e., $\pi _{1}\left( 1\right) =1$
and $\pi _{1}\left( 4\right) =4$.

Now lets consider the blocks $A_{2}$ and $A_{4},$ i.e.$,$the equivalent
classes $[2,3]$ and $[5,6],$ and consider the automorphism $\pi _{1}$ such
that $\pi _{1}\left( 3\right) =6$ and $\pi _{1}\left( 4\right) =4$ then $%
\left( \pi _{1}\left( 3\right) ,\pi _{1}\left( 4\right) \right) =\left(
6,4\right) \notin \rho ^{\prime }.$ Therefore these blocks can not permute
so the equivalent classes are not comparable, i.e., $\pi _{1}\left(
[2,3]\right) =[2,3]$ and $\pi _{1}\left( [5,6]\right) =[5,6]$. So $P_{\pi
_{1}}=I.$

Consider the permutation $\pi _{1}=(1,4)\in $ $\rho ^{\prime }$ and non
trivial transitive function

$g(i;j)=\left\{ 
\begin{array}{c}
a\text{ if }\left( i,j\right) =\left( 1,4\right) \\ 
1\text{ if }\left( i,j\right) \neq \left( 1,4\right)%
\end{array}%
\right. $ where $a\in F$ and $a\neq 1$.

Then we have an automorphism $G:M_{6}\left( F,\rho ^{\prime }\right)
\rightarrow $ $M_{6}\left( F,\rho ^{\prime }\right) $ generated by
transitive function $g(i;j)$ as follows

$G(i;j)=\left\{ 
\begin{array}{c}
aE^{14}\text{ if }\left( i,j\right) =\left( 1,4\right) \\ 
E^{ij}\text{ \ if }\left( i,j\right) \neq \left( 1,4\right)%
\end{array}%
\right. .$

We have an equality of $G$ and $\Gamma $ which is obtained in the proof of
the Factorization Theorem (see \cite{Akkurt E-M Barker 2}) that is%
\begin{equation*}
\Gamma \left( E^{ij}\right) =G\left( E^{\pi _{1}\left( i\right) \pi
_{1}\left( j\right) }\right) =\left( G\circ P_{\pi _{1}}\right) \left(
E^{ij}\right) ,
\end{equation*}%
and since $P_{\pi _{1}}=1$ then we have $\Gamma =G.$

So

\begin{equation*}
\Phi =\Psi _{\overline{A}}\circ G.
\end{equation*}

\end{document}